\input amstex
\documentstyle{amsppt}
\magnification1200
\def\phi{\varphi}
\def\0{\varnothing}
\def\eps{\varepsilon}

\def\R{\Bbb R}
\def\N{\Bbb N}

\def\bs{\backslash}

\topmatter
\title
On cake dividing
\endtitle
\rightheadtext{On cake dividing}
\author
Olexandr Ravsky  
\endauthor
\email oravsky\@mail.ru
\endemail
\abstract
Considered Steinhaus geometric game on cake dividing.
\endabstract
\endtopmatter
\document
 Hugo Steinhaus in his popular book [S] considered the following game
(problem 51). Pavel and Havel are dividing a cake as follows. Pavel
selects a point $P$ and Havel draws a line $l$ through the point $P$ and
gets his piece of the cake. What form of the cake is most advantageous for
Havel and which part of the cake he obtains in this case? 

Now we consider the formal interpretation of the game.  

Let $\mu$ be the Lebesgue measure on the space $\R^n$.
Let $\Cal B^n$ be a family of compacts of the space $\R^n$ with non zero 
measure. For every point $a\in\R^n$ put
$l_a=\{x\in\R^n:(x,a)\ge 0\}$. Put 
$$b_n=\inf_{A\in\Cal B^n}\sup_{x\in A}\inf_{a\in\R^n}
\frac{\mu(A\cap(x+l_a))}{\mu(A)}.$$

Steinhaus proved that $1/4\le b_2\le 1/3$. 

\proclaim{Theorem 1} For every $n$ holds $b_n=1/(n+1)$.
\endproclaim
\demo{Proof} Let $S\subset\R^n$ be the regular simplex with vertices 
$a_0,\dots a_n$. Every point $a\in S$ has an unique representation 
$a=\sum \lambda_i a_i$ such that $\lambda_i\ge 0$ for every $i$ and 
$\sum\lambda_i=1$. Put $S_i=a_i-S$ for every $i$. If $i\not=j$ then
$S_i\cap S_j=\{0\}$. Indeed suppose that there exist points $s_i,s_j\in S$
such that $x=a_i-s_i=a_j-s_j$. Let $s_j=\sum\lambda_{ik} a_k$ and 
$s_j=\sum\lambda_{jk} a_k$ be the representations of the points $s_i$ and
$s_j$. The uniquity of the representation of the point
$(a_i+s_j)/2=(a_j+s_i)/2$ gives that $1+\lambda_{ji}=\lambda_{ii}$.
Thus $\lambda_{ii}=1$, $s_i=a_i$ and $x=0$. 

Put $A=\bigcup S_i$. Then $\mu(A)=(n+1)\mu(S)$. Let $a\in S$,
$a=\sum \lambda_i a_i$ be the representation of the point $a$ and
$\lambda_j=\max\lambda_i$. We show that $S_j\subset l_a$. Indeed let
$a_j-x=\sum \mu_i a_i$ be the representation of a point $a_j-x$ where
$x\in S_j$. Let $\alpha=(a_i,a_j)$ where $i\not=j$. 
Then $(x,a_j)=1-\mu_j-\alpha\sum_{i\not=j}\mu_i\ge 1-\sum\mu_i\ge 0.$
Thus the construsted set $A$ shows that $b_n\le 1/(n+1)$.   

Show now that $b_n\ge 1/(n+1)$. Consider an arbitrary set $A\in\Cal B^n$.
For every $a\in S$ and $0\le\eps\le 1/(n+1)$ put 
$A(a,\eps)=\{x\in\R^n:\mu(A\cap(x+l_a))\ge 1/(n+1)-\eps\}$. 
Then every $A(a,\eps)$ is a convex closet set and $\mu(A(a,\eps)\cap A)=n/(n+1)+\eps$.
Now fix an arbitrary number $\eps>0$. Let $c_0,\dots,c_n\in S$. 
Since $\mu(A)>\sum\mu(A\bs A(c_i,\eps))$ then there exists a point 
$c(\eps)\in A\cap\bigcap A(c_i,\eps)$.
Since $A$ is a compact then there exists a cluster point $c$ of 
the set $\{c(1/m):m\in\N\}$. Then $c\in A\cap\bigcap A(c_i,0)$.
Then the Helly Theorem implies that $\bigcap\{A(c,0):c\in\R^n\}\not=\0$. 
Therefore $b_n\ge 1/(n+1)$.
$\square$
\enddemo

\heading
References
\endheading

[G] B. Gr\"unbaum, {\it Etudes on combinatorial geometry and theory of convex bodies},
M., Nauka, 1971. (in Russian)

[S] H. Steinhaus, {\it Sto zada\'n}, M., Nauka, 1982. (in Russian)
\enddocument
\bye